\documentclass[11pt]{article}
\usepackage{amssymb,amsmath,latexsym, verbatim}
\usepackage{tikz}
\allowdisplaybreaks

\begin{document}

\begin{doublespace}

\newtheorem{thm}{Theorem}[section]
\newtheorem{theorem}{Theorem}[section]
\newtheorem{lemma}[thm]{Lemma}
\newtheorem{cond}[thm]{Condition}
\newtheorem{defn}[thm]{Definition}
\newtheorem{prop}[thm]{Proposition}
\newtheorem{proposition}[thm]{Proposition}
\newtheorem{corollary}[thm]{Corollary}
\newtheorem{remark}[thm]{Remark}
\newtheorem{example}[thm]{Example}
\newtheorem{conj}[thm]{Conjecture}
\numberwithin{equation}{section}
\def\ee{\varepsilon}
\def\qed{{\hfill $\Box$ \bigskip}}
\def\NN{{\cal N}}
\def\AA{{\cal A}}
\def\MM{{\cal M}}
\def\BB{{\cal B}}
\def\CC{{\cal C}}
\def\LL{{\cal L}}
\def\DD{{\cal D}}
\def\FF{{\cal F}}
\def\EE{{\cal E}}
\def\QQ{{\cal Q}}
\def\RR{{\mathbb R}}
\def\R{{\mathbb R}}
\def\L{{\bf L}}
\def\K{{\bf K}}
\def\S{{\bf S}}
\def\A{{\bf A}}
\def\E{{\mathbb E}}
\def\F{{\bf F}}
\def\P{{\mathbb P}}
\def\Prob{{\mathbb P}}

\def\N{{\mathbb N}}
\def\eps{\varepsilon}
\def\wh{\widehat}
\def\wt{\widetilde}
\def\pf{\noindent{\bf Proof.} }
\def\beq{\begin{equation}}
\def\eeq{\end{equation}}
\def\lam{\lambda}
\def\H{\mathcal{H}}
\def\nn{\nonumber}
\def\C{\mathbb{C}}

\newcommand\blfootnote[1]{%
  \begingroup
  \renewcommand\thefootnote{}\footnote{#1}%
  \addtocounter{footnote}{-1}%
  \endgroup
}

\title{\Large \bf  Spectral heat content for  $\alpha$-stable processes in $C^{1,1}$ open sets}
\author{Hyunchul Park\thanks{Research  supported in part by Research and Creative Project Award from SUNY New Paltz.}  and Renming Song\thanks{Resarch supported in part by a grant from the Simons Foundation (\#429343, Renming Song).}}

\date{ }
\maketitle

\begin{abstract}
In this paper we study the asymptotic behavior, as $t\downarrow 0$, of the spectral heat content 
$Q^{(\alpha)}_{D}(t)$ for 
isotropic $\alpha$-stable processes, $\alpha\in [1,2)$, in bounded $C^{1,1}$ open sets 
$D\subset \R^{d}$, $d\geq 2$.
Together with the results from \cite{Val2017} for $d=1$ and \cite{GPS19} for $\alpha\in (0,1)$,
the main theorem of this paper establishes the asymptotic behavior of the spectral heat content 
up to the second term for all $\alpha\in (0,2)$ and $d\geq1$, 
and resolves the conjecture raised in \cite{Val2017}. 
\end{abstract}

\section{Introduction}
The spectral heat content represents the total heat 
in a domain $D$ with Dirichlet boundary condition when the initial temperature is 1.
The spectral heat content for Brownian motions has been studied extensively.
The spectral heat content for isotropic stable processes was first studied in \cite{Val2017}. 
Since then, considerable progress has been made toward understanding the asymptotic behavior of the spectral heat content for other L\'evy processes (see \cite{Val2016,  GPS19, P20, PS19}).

The following conjecture about the spectral heat content for isotropic $\alpha$-stable
processes, $\alpha\in (0,2)$, over bounded $C^{1,1}$ open sets (see Section \ref{section:preliminaries} for the definition of $C^{1,1}$ open sets) was made in \cite{Val2017}: 
As $t\downarrow 0$,
\beq\label{eqn:Val conj}
Q_{D}^{(\alpha)}(t)=
\begin{cases}
|D|-c_{1}|\partial D|t^{1/\alpha} +O(t), \quad \alpha\in (1,2),\\
|D|-c_{2}|\partial D|t\ln(1/t) +O(t), \quad \alpha=1,\\
|D|-c_{3}\text{Per}_{\alpha}(D)t +o(t), \quad \alpha\in (0,1),
\end{cases}
\eeq
where $c_i, i=1, 2, 3,$ are constants and 
\beq\label{eqn:fractional perimeter}
\text{Per}_{\alpha}(D):=\int_{D}\int_{D^{c}}\frac{A_{d,\alpha}}{|x-y|^{d+\alpha}}dydx,
\quad \mbox{ with }A_{d,\alpha}=\frac{\alpha\Gamma(\frac{d+\alpha}{2})}{2^{1-\alpha}\pi^{d/2}\Gamma(1-\frac{\alpha}{2})},
\eeq
is the $\alpha$-fractional perimeter of $D$.
This conjecture was resolved in dimension 1 in \cite{Val2016} (actually, a slightly weaker version, with the error term being $o(t^{1/\alpha})$ in the case $\alpha\in (1, 2)$ and $o(t\ln(1/t))$ in the case $\alpha=1$, of the conjecture was proved there).
We note that in \cite{Val2017} the author 
conjectured, and also provided strong evidence, that the spectral heat content for isotropic $\alpha$-stable processes  with $\alpha\in (0,2)$ must have an asymptotic expansion of the form as \eqref{eqn:Val conj} for all dimensions $d\geq2$, 
but exact expressions for the coefficients $c_{i}$ were not provided.
Then a two-term asymptotic expansion of the spectral heat content for L\'evy processes of bounded variation in $\R^{d}$ was established in \cite{GPS19}.
Since $\alpha$-stable processes are of bounded variation if and only if $\alpha\in (0,1)$, 
the result in \cite{GPS19} proves \eqref{eqn:Val conj} for $\alpha\in (0,1)$. 
The purpose of this paper is to resolve the conjecture above for $\alpha \in [1,2)$ and $d\geq 2$. 
In fact, our result is slightly weaker than \eqref{eqn:Val conj} since the error term is $o(t^{1/\alpha})$ for $\alpha\in (1, 2)$ and $o(t\ln(1/t))$ for $\alpha=1$.
We also find explicit expressions for the constants $c_{1}$ and $c_{2}$.
The main results of this paper are Theorem \ref{thm:main12} for $\alpha\in (1,2)$ and Theorem \ref{thm:main1} for $\alpha=1$. 
Combining Theorems \ref{thm:main12} and \ref{thm:main1} with \cite[Corollary 3.5]{GPS19}, the asymptotic behavior of the spectral heat content for isotropic $\alpha$-stable processes in bounded $C^{1,1}$ open sets $D$ can be stated as follows: 

\begin{theorem}\label{thm:main}
Let $D$ be a bounded $C^{1,1}$ open set in $\R^{d}$, $d\geq 2$ and 
let
$$
f_{\alpha}(t)=
\begin{cases}
t^{1/\alpha} &\text{ if } \alpha\in (1,2),\\
t\ln(1/t) &\text{ if } \alpha=1,\\
t &\text{ if } \alpha\in (0,1).
\end{cases}
$$
Then, we have
\beq\label{eqn:main result}
\lim_{t\to 0}\frac{|D|-Q_{D}^{(\alpha)}(t)}{f_{\alpha}(t)}=
\begin{cases}
\E[\overline{Y}_{1}^{(\alpha)}]|\partial D| &\text{if } \alpha\in (1,2),\\
\frac{|\partial D|}{\pi} &\text{if }\alpha=1,\\
\rm{Per}_{\alpha}(D), &\text{if } 
\alpha\in (0,1),
\end{cases}
\eeq
where $\overline{Y}_{t}^{(\alpha)}=\sup_{s\le t}Y^{(\alpha)}_s$ stands for the running supremum of a 1 dimensional symmetric $\alpha$-stable process $Y^{(\alpha)}_t$ and 
$\rm{Per}_{\alpha}(D)$ as defined in \eqref{eqn:fractional perimeter}.
\end{theorem}
We note that \eqref{eqn:main result} is exactly the same form as \cite[Theorem 1.1]{Val2016} if one interprets $|\partial D|=2$ when $D$ is a bounded open interval in $\R$, but the proof for $d\geq 2$ is very different from the one dimensional case and much more challenging.

The two-term asymptotic expansion of the spectral heat content for Brownian motion was proved in \cite{BD89}. 
The crucial ingredient in \cite{BD89} is the fact that 
individual components of Brownian motion are independent.
For isotropic $\alpha$-stable processes, $\alpha\in (0, 2)$, individual components  are not independent and the technique in \cite{BD89} no longer works.
 When $\alpha\in (1,2)$,
we establish the lower bound for the heat loss $|D|-Q_{D}^{(\alpha)}(t)$ by considering the most efficient way of exiting $D$ (see Lemma \ref{lemma:lb}).
In order to establish the upper bound, we  approximate the heat loss 
$|D|-Q_{D}^{(\alpha)}(t)$ by the heat loss of the half-space 
$|D|-Q_{H_{x}}^{(\alpha)}(t)$ for  
$x$ near $\partial D$ (see Proposition \ref{prop:ub12}) 
and show that the approximation error is of order $o(t^{1/\alpha})$ in Lemma \ref{lemma:near boundary}, 
which is similar to tools exposed in the trace estimate results in \cite{BK08, PS14}.
However, these tools do not work when $\alpha=1$ due to the non-integrability of 
$\Prob(\overline{Y}^{(1)}_{1}>u)$ over $(1,\infty)$, where $ \overline{Y}^{(1)}_t$ stands for the
supremum of the Cauchy process up to time $t$, 
and the proof for $\alpha=1$ requires new ideas and is considerably more difficult. 

In case of $\alpha=1$, we prove that the coefficient of the second term 
of the asymptotic expansion of $Q_{D}^{(1)}(t)$ is $-\frac{|\partial D|}{\pi}$, 
which is the same as for the regular heat content $H^{(1)}_{D}(t)$ 
which represents the total heat in $D$ without the Dirichlet exterior condition 
(see \eqref{e:rhc} below for the definition of regular heat content and \cite[Theorem 1.2]{Val2017} 
for the two-term asymptotic expansion for the regular heat content $H_{D}^{(1)}(t)$). 
If there is no Dirichlet exterior condition on $D^{c}$, or equivalently if the heat moves freely in and out of $D$,
the heat loss of the regular heat content must be smaller than that of the spectral heat content and we obtain the lower bound for free in \eqref{eqn:Cauchy lb}. The proof for the upper bound is much more demanding. 
The crucial ingredient for the upper bound is the spectral heat content for subordinate killed Brownian motions in \cite{PS19}. 
An isotropic $\alpha$-stable process $X_{t}^{(\alpha)}$ can be realized as a subordinate Brownian motion $W_{S_{t}^{(\alpha/2)}}$, where $S_{t}^{(\alpha/2)}$ is an independent $(\alpha/2)$-stable subordinator.
Hence, the spectral heat content $Q_{D}^{(\alpha)}(t)$ is the spectral heat content for the \textit{killed subordinate Brownian motion} via the independent $(\alpha/2)$-stable subordinator $S_{t}^{(\alpha/2)}$. 
When one reverses the order of killing and subordination, one obtains the \textit{subordinate killed Brownian motion}. 
This is the process obtained by subordinating the killed Brownian motion in $D$ via the independent $(\alpha/2)$-stable subordinator $S_{t}^{(\alpha/2)}$. 
Let $\widetilde{Q}_{D}^{(\alpha)}(t)$ be the spectral heat content of the subordinate killed Brownian motion in $D$
(see \eqref{eqn:SKBM} for the precise definition).
By construction,
 $\widetilde{Q}_{D}^{(\alpha)}(t)$ is always smaller than $Q_{D}^{(\alpha)}(t)$, 
or equivalently the heat loss $|D|-\widetilde{Q}_{D}^{(\alpha)}(t)$ of the subordinate killed Brownian motion provides a natural upper bound for 
that of the killed subordinate Brownian motion $|D|-Q_{D}^{(\alpha)}(t)$. 
We use two independent $\alpha/2$ and $\beta/2$ stable subordinators 
with $\alpha\beta=2$ 
and consider the $\alpha$-stable process $X^{(\alpha)}_{t}=W_{S_{t}^{(\alpha/2)}}$ and killed it upon exiting $D$. 
Then we time-change the killed $\alpha$-stable process by the $(\beta/2)$-stable subordinator. 
By using the heat loss of the resulting process, we obtain in Lemma 
\ref{lemma:Cauchy ub} that, for $\alpha\in (1,2)$, 
$$
\limsup_{t\to 0}\frac{|D|-Q_{D}^{(1)}(t)}{t\ln(1/t)}\leq
\frac{|\partial D|}{\pi} +
\frac{\int_{0}^{\infty}\Prob(\overline{Y}^{(\alpha)}_{1}>u, Y_{1}^{(\alpha)}<u )du }
{\Gamma(1-\frac{1}{\alpha})}\cdot |\partial D|,
$$
where $Y^{(\alpha)}_t$ is a 1 dimensional symmetric $\alpha$-stable process
and $\overline{Y}_{t}^{(\alpha)}=\sup_{s\le t}Y^{(\alpha)}_s$. 
We would like to show that the second summand on the right hand side of the previous inequality converges to 0 as $\alpha\downarrow 1$. 
We know that
$\frac{1}{\Gamma(1-\frac{1}{\alpha})}=O(\alpha-1)$ as $\alpha\downarrow 1$. 
The integrand in the  numerator can be written as 
$$
\Prob(\overline{Y}^{(\alpha)}_{1}>u, Y_{1}^{(\alpha)}<u )du
=\Prob(\overline{Y}^{(\alpha)}_{1}>u)-\Prob(Y^{(\alpha)}_{1}\geq u)
$$
and it can be shown that 
$\Prob(\overline{Y}^{(\alpha)}_{1}>u)\sim \Prob(Y^{(\alpha)}_{1}\geq u) \sim cu^{-\alpha}$ 
as $u\to \infty$, 
hence
$\int_{0}^{\infty}\Prob(\overline{Y}^{(\alpha)}_{1}>u)du$ and $\int_{0}^{\infty}\Prob(Y^{(\alpha)}_{1}\geq u)du$ 
should be of order $\frac{1}{\alpha-1}$ as $\alpha\downarrow 1$.
We show that
$\int_{0}^{\infty}\Prob(\overline{Y}^{(\alpha)}_{1}>u)du$ and $\int_{0}^{\infty}\Prob(Y^{(\alpha)}_{1}\geq u)du$ 
have \textit{the exactly same leading coefficients} and this gives a cancellation of the main terms of order $\frac{1}{\alpha-1}$. 
We still need to show that the sub-leading terms are of order $o(\alpha-1)$ as $\alpha\downarrow 1$ and we show this by establishing uniform heat kernel upper and lower bounds in Lemma \ref{lemma:uniform ub} 
for the heat kernel  of the supremum process $\overline{Y}_{t}^{(\alpha)}$ 
and Lemma \ref{lemma:uniform lb} for the heat kernel of $Y_{t}^{(\alpha)}$.

This paper deals with asymptotic behavior of the spectral heat content for isotropic $\alpha$-stable processes. 
It is natural and interesting to try to find the  asymptotic behavior of the spectral heat content for more general  L\'evy processes. We intend to deal with this topic in a future project.
In the recent paper \cite{P20}, a three-term asymptotic expansion of the spectral heat content of 1 dimensional symmetric $\alpha$-stable processes, $\alpha\in [1,2)$, was established.
We believe that a similar result should hold true for $d\geq 2$.

The organization of this paper is as follows. 
In Section \ref{section:preliminaries}, we introduce the setup. In Section \ref{section:alpha12},
we deal with the case $\alpha\in (1,2)$ and the main result of that section is Theorem \ref{thm:main12}.  
The case $\alpha=1$ is dealt with in Section \ref{section:alpha1} and  the main result there is Theorem \ref{thm:main1}. 
In this paper, we use $c_i$ to denote constants whose values are unimportant and may change from one appearance to another.
The notation $\Prob_{x}$ stands for the law of the underlying processes started at $x\in \R$, and $\E_{x}$ stands for expectation with respect to $\Prob_{x}$. For simplicity, we use $\Prob=\Prob_{0}$ and $\E=\E_{0}$.

\section{Preliminaries}\label{section:preliminaries}

In this paper,  unless explicitly stated otherwise, we assume $d\geq 2$.
Let $X_{t}^{(\alpha)}$, $\alpha\in (0, 2]$, be an isotropic $\alpha$-stable process with 
$$
\E[e^{i\xi X^{(\alpha)}_{t}}]=e^{-t|\xi|^{\alpha}}, \quad \xi \in \R^{d}, \alpha\in (0,2].
$$
$X_{t}^{(2)}$ is a Brownian motion $W_t$ with transition density given by $(4\pi t)^{-d/2}e^{-\frac{|x|^2}{4t}}$. 
Let $S_{t}^{(\alpha/2)}$, $\alpha\in (0,2)$, be an $(\alpha/2)$-stable subordinator with
$$
\E[e^{-\lam S_{t}^{(\alpha/2)}}]=e^{-t \lam^{\alpha/2}}, \quad \lam>0.
$$
Assume that $S_{t}^{(\alpha/2)}$ is independent of the Brownian motion $W_t$. 
Then, the subordinate Brownian motion $W_{S_{t}^{(\alpha/2)}}$ is a realization of the process $X_{t}^{(\alpha)}$.
We will reserve $Y_{t}^{(\alpha)}$ for the 1 dimensional symmetric $\alpha$-stable process.
We define the running supremum process $\overline{Y}^{(\alpha)}_{t}$ of  $Y_{t}^{(\alpha)}$ by
\beq\label{eqn:running}
\overline{Y}^{(\alpha)}_{t}=\sup\{Y^{(\alpha)}_{u} : 0\leq u\leq t\}.
\eeq

Recall that an open set $D$ in $\R^d$ is said to be a $C^{1, 1}$ open set if 
there exist a localization radius $R_0>0$ and a constant $\Lambda_0>0$ such that, 
for every $z\in \partial D$, there exist a $C^{1, 1}$ function $\phi=\phi_z: \R^{d-1}\to \R$ satisfying $\phi(0)=0$, $\nabla \phi(0)=(0, \cdots, 0)$, $\|\nabla\phi\|_\infty\le \Lambda_0$,
$|\nabla \phi(x)-\nabla \phi(y)|\le \Lambda_0 |x-y|$ and an orthonormal coordinate system $CS_z: y=(\widetilde y, y_d)$
with origin at $z$ such that
$$
B(z, R_0)\cap D=B(z, R_0)\cap \{ y=(\widetilde y, y_d) \mbox{ in } CS_z: y_d>\phi(\widetilde y)\}.
$$
The pair $(R_0, \Lambda_0)$ is called the $C^{1,1}$ characteristics of the $C^{1, 1}$ open set $D$. 
It is well known that any $C^{1, 1}$ open set $D$  in $\R^d$ satisfies the uniform interior and exterior $R$-ball condition: for any $z\in\partial D$,  there exist balls $B_1$ and $B_2$ of radii $R$ with $B_1\subset D$, $B_2\subset\R^d\setminus\overline{D}$ and $\partial B_1\cap\partial B_2=\{z\}$.

We recall from \cite{BD89} a useful fact about open sets $D$ satisfying the uniform interior and exterior $R$-ball condition.
Let $D_{q}=\{x\in D : \text{dist}(x, \partial D) >q\}$.
We will use $\partial D_{q}$ denote the portion of the boundary of $D_q$ contained in $D$, that is, $\partial D_{q}=\{x\in D: \text{dist}(x, \partial D)=q\}$.
It follows from \cite[Lemma 6.7]{BD89}
that
\beq\label{eqn:stability}
|\partial D|\left(\frac{R-q}{R}\right)^{d-1}\leq |\partial D_{q}|\leq |\partial D|\left(\frac{R}{R-q}\right)^{d-1}, \quad 0< q<R.
\eeq

In the remainder of this paper, $D$ stands for a bounded $C^{1,1}$ open set in $\R^{d}$, $d\geq 2$.
Let $\tau_{D}^{(\alpha)}=\inf\{t>0: X^{(\alpha)}_{t}\notin D\}$ be the first exit time of $X^{(\alpha)}_{t}$ from $D$. 
The spectral heat content  of $D$ for $X^{(\alpha)}_{t}$  is defined to be
$$
Q_{D}^{(\alpha)}(t)=\int_{D}\Prob_{x}(\tau_{D}^{(\alpha)}>t)dx.
$$
The (regular) heat content $H_{D}^{(\alpha)}(t)$ of $D$  for $X^{(\alpha)}_{t}$  is defined to be
\begin{equation}\label{e:rhc}
H_{D}^{(\alpha)}(t)=\int_{D}\Prob_{x}(X_{t}^{(\alpha)}\in D)dx.
\end{equation}
The spectral heat content $\widetilde{Q}_{D}^{(\alpha)}(t)$  of $D$ for the subordinate killed Brownian motion is defined as 
\beq\label{eqn:SKBM}
\widetilde{Q}_{D}^{(\alpha)}(t)=\int_{D}\Prob_{x}(\tau_{D}^{(2)}>S_{t}^{(\alpha/2)})dx.
\eeq
Note that we always have 
\beq\label{eqn:two heat contents}
\{\tau_{D}^{(2)}>S_{t}^{(\alpha/2)}\} \subset \{\tau_{D}^{(\alpha)}>t\}\subset \{X_{t}^{(\alpha)}\in D\}.
\eeq
We sometimes use the terminology \textit{heat loss} and this will mean either
$$
|D|-Q_{D}^{(\alpha)}(t)=\int_{D}\Prob_{x}(\tau_{D}^{(\alpha)}\leq t)dx,
$$
or 
$$
|D|-\widetilde{Q}_{D}^{(\alpha)}(t)=\int_{D}\Prob_{x}(\tau_{D}^{(2)}\leq S_{t}^{(\alpha/2)})dx,
$$
depending on which process we are dealing with. 
Intuitively, these quantities represent the total heat loss caused by heat particles jumping out of $D$ up to time $t$.
From \eqref{eqn:two heat contents}, we have 
$$
|D|- H_{D}^{(\alpha)}(t)\leq |D|-Q_{D}^{(\alpha)}(t) \leq |D|-\widetilde{Q}_{D}^{(\alpha)}(t), \quad t>0.
$$

\section{The case $\alpha\in (1,2)$}\label{section:alpha12}
Throughout this section, we assume $\alpha\in(1,2)$. The main purpose of this section is to prove Theorem \ref{thm:main12}.

For any $x\in D$, we use $\delta_D(x)$ to denote the distance between $x$ and $\partial D$.
For a bounded $C^{1,1}$ open set $D$ satisfying the uniform 
interior and exterior $R$-ball condition
and $x\in D$ with $\delta_{D}(x)<R/2$, we let $z_x\in \partial D$ be the point on $\partial D$ such that $|x-z_x|=\delta_D(x)$ and 
$H_{x}$ be the half-space containing the interior 
$R$-ball at the point  $z_x$ and tangent to $\partial D$ at $z_x$.
For $x\in D$ with $\delta_{D}(x)<R/2$, we have
\begin{eqnarray}\label{eqn:ub cases2}
\{\tau_{D}^{(\alpha)}\leq t\} 
&=& \{\tau_{D}^{(\alpha)}\leq t \text{ and } \tau_{H_{x}}^{(\alpha)}\leq t\} \cup \{\tau_{D}^{(\alpha)}\leq t\text{ and } \tau_{H_{x}}^{(\alpha)}>t\}\nn\\
&\subset&\{\tau_{H_{x}}^{(\alpha)}\leq t\} \cup \{\tau_{D}^{(\alpha)}\leq t <\tau_{H_{x}}^{(\alpha)}\}.
\end{eqnarray}
Hence, we have
\begin{eqnarray}\label{eqn:ub cases}
&&|D|-Q^{(\alpha)}_{D}(t)=\int_{D}\Prob_{x}(\tau^{(\alpha)}_{D}\leq t)dx=\int_{D_{R/2}}\Prob_{x}(\tau_{D}^{(\alpha)}\leq t)dx +\int_{D\setminus D_{R/2}}\Prob_{x}(\tau_{D}^{(\alpha)}\leq t)dx\nn\\
&\leq &\int_{D_{R/2}}\Prob_{x}(\tau^{(\alpha)}_{D}\leq t)dx+\int_{D\setminus D_{R/2}}\Prob_{x}(\tau^{(\alpha)}_{H_{x}}\leq t)dx +\int_{D\setminus D_{R/2}}\Prob_{x}(\tau^{(\alpha)}_{D}\leq t<\tau^{(\alpha)}_{H_{x}})dx.
\end{eqnarray}

We deal with the first expression of \eqref{eqn:ub cases} first.
Recall the following facts from \cite[(2.2)]{Val2017} and \cite[Corollary 6.4]{BD89}:
\beq\label{eqn:exp moment}
\E\left[\exp\left(-\frac{\kappa^{2}}{S_{1}^{(\alpha/2)}}\right)\right]\leq c_{\alpha}\kappa^{-\alpha}
\eeq
and
\beq\label{eqn:BM exit}
\Prob\left(\tau^{(2)}_{B(0,R)}<t\right)\leq 2^{1+d/2}e^{-\frac{R^2}{8t}}.
\eeq

\begin{lemma}\label{lemma:inside}
Let $D\subset \R^{d}$, $d\geq 2$, be a bounded $C^{1,1}$ open set satisfying the uniform 
interior and exterior $R$-ball condition. 
Then, we have
$$
\int_{D_{R/2}}\Prob_{x}\left(\tau_{D}^{(\alpha)}\leq t\right)dx\leq ct.
$$
\end{lemma}
\pf
For $\delta_{D}(x)\geq \frac{R}{2}$, it follows from \eqref{eqn:exp moment} and \eqref{eqn:BM exit} that 
\begin{align*}
&\Prob_{x}\left(\tau_{D}^{(\alpha)}\leq t\right)
\leq \Prob_{x}\left(\tau^{(\alpha)}_{B(x,\delta_{D}(x))} \leq t\right)
\leq \Prob_{x}\left(\tau^{(2)}_{B(x,\delta_{D}(x))} \leq S^{(\alpha/2)}_{t}\right)\\
\leq& 2^{1+\frac{d}2}\E[\exp(-\frac{\delta_{D}(x)^{2}}{8S_{t}^{(\alpha/2)}})]
\leq 2^{1+\frac{d}2}\E\left[\exp(-\frac{\delta_{D}(x)^{2}}{8t^{2/\alpha}S_{1}^{(\alpha/2)}})\right]
\leq\frac{ c(\alpha,d)t}{\delta_{D}(x)^{\alpha}}
\leq \frac{2^{\alpha}c(\alpha,d)t}{R^{\alpha}}.
\end{align*}
Hence we have
$$
\int_{D_{R/2}}\Prob_{x}\left(\tau_{D}^{(\alpha)}\leq t\right)dx\leq \frac{2^{\alpha}c(\alpha,d)|D|t}{R^{\alpha}}.
$$
\qed

The second expression of \eqref{eqn:ub cases} is handled in the following proposition. 
\begin{proposition}\label{prop:ub12}
Let $D\subset \R^{d}$, $d\geq 2$, be a bounded $C^{1,1}$ open set satisfying the uniform 
interior and exterior $R$-ball condition. 
Then, we have
$$
\lim_{t\to 0}\frac{\int_{D\setminus D_{R/2}}\Prob_{x}(\tau^{(\alpha)}_{H_{x}}\leq t)dx}{t^{1/\alpha}}=|\partial D|
\int_{0}^{\infty}\Prob_{(0,r)}(\tau^{(\alpha)}_{H}\leq 1)dr=|\partial D|
\E[\overline{Y}^{(\alpha)}_{1}],
$$
where $\overline{Y}_{t}^{(\alpha)}$ is defined in \eqref{eqn:running} and $H=\{x=(x_{1},\cdots ,x_{d}) : x_{d}>0\}$.
\end{proposition}
\pf
For $\delta_{D}(x)<R/2$, we have by the scaling property and the change of variables $v=t^{-1/\alpha}u$,
\begin{align*}
&\int_{D\setminus D_{R/2}}\Prob_{x}\left(\tau^{(\alpha)}_{H_{x}}\leq t\right)dx
=\int_{0}^{R/2}|\partial D_{u}|
\Prob_{(\widetilde 0, u)}\left(\tau_{H}^{(\alpha)}\leq t\right)du\\
&=\int_{0}^{R/2}|\partial D_{u}|
\Prob_{(\widetilde 0, t^{-1/\alpha}u)}
\left(\tau_{t^{-1/\alpha}H}^{(\alpha)}\leq 1\right)du
=t^{1/\alpha}\int_{0}^{R/2t^{1/\alpha}}|\partial D_{vt^{1/\alpha}}|
\Prob_{(\widetilde 0, v)}\left(\tau_{H}^{(\alpha)}\leq 1\right)dv.
\end{align*}

Recall that $Y^{(\alpha)}_t=(X_{t}^{(\alpha)}-x)\cdot {\bf n}_{z_x}$ is a one dimensional $\alpha$-stable process starting from 0. 
Hence, for any starting point $(\tilde{0},r)$  
we have $\{\tau^{(\alpha)}_{H}\leq 1\}=\{\overline{Y}^{(\alpha)}_{1}>r\}$, 
and this implies
$\Prob_{(\widetilde 0, r)}(\tau^{(\alpha)}_{H}\leq 1)=
\Prob(\overline{Y}^{(\alpha)}_{1}>r)$. 
Hence, we have
$$
\int_{0}^{\infty}\Prob_{(\widetilde 0, v)}(\tau^{(\alpha)}\leq 1)dv
=\int_{0}^{\infty}\Prob(
\overline{Y}^{(\alpha)}_{1}\geq v)dv=\E[\overline{Y}^{(\alpha)}_{1}].
$$
By \cite[Corollary 2.1 (ii)]{Val2016}, $\E[\overline{Y}^{(\alpha)}_{1}]<\infty$ if $\alpha\in (1,2)$.
Hence, it follows from the Lebesgue dominated convergence theorem that
$$
\lim_{t\to 0}\frac{\int_{D\setminus D_{R/2}} \Prob_{x}(\tau^{(\alpha)}_{H_{x}} \leq t)dx}{t^{1/\alpha}}=|\partial D|\int_{0}^{\infty}
\Prob_{(\widetilde 0, r)}(\tau^{(\alpha)}_{H}\leq 1)dr=
|\partial D|\E[\overline{Y}^{(\alpha)}_{1}].
$$
\qed

Finally, we estimate the last expression of \eqref{eqn:ub cases}.
\begin{lemma}\label{lemma:near boundary}
Let $D\subset \R^{d}$, $d\geq 2$, be a bounded $C^{1,1}$ 
open set. 
Then, we have
$$
\lim_{t\to 0}\frac{\int_{D\setminus D_{R/2}}\Prob_{x}\left(\tau_{D}^{(\alpha)}\leq t< \tau^{(\alpha)}_{H_{x}}\right)dx}{t^{1/\alpha}}=0.
$$
\end{lemma}
\pf
By rotational invariance, we have
\begin{align}\label{eqn:approx1}
&\int_{D\setminus D_{R/2}}\Prob_{x}\left(\tau_{D}^{(\alpha)}\leq t< \tau^{(\alpha)}_{H_{x}}\right)dx
\leq\int_{D\setminus D_{R/2}}\Prob_{x}\left(\tau^{(\alpha)}_{B_{x}}\leq t< \tau^{(\alpha)}_{H_{x}}\right)dx\nn\\
&\leq \int_{0}^{R/2}|\partial D_{u}|\Prob_{(\widetilde 0, u)}
\left(\tau^{(\alpha)}_{B((\widetilde 0, R), R)}
\leq t< \tau^{(\alpha)}_{H}\right)du,
\end{align}
where $H=\{x=(x_{1},\cdots ,x_{d}) : x_{d}>0\}$.
It follows from \cite[Lemma 6.7]{BD89} that  for $u<R/2$,
\beq\label{eqn:approx2}
|\partial D_{u}| \leq 2^{d-1}|\partial D|.
\eeq

Note that it follows from the scaling property that the law of $a\tau^{(\alpha)}_{D}$ under $\Prob_{x}$ is equal to the law of $\tau^{(\alpha)}_{a^{1/\alpha}D}$ under $\Prob_{a^{1/\alpha}x}$.
Hence, it follows from the scaling property, together with  \eqref{eqn:approx2}, that \eqref{eqn:approx1} is bounded above by
\begin{eqnarray}\label{eqn:approx3}
&&2^{d-1}|\partial D|\int_{0}^{R/2}
\Prob_{(\widetilde 0,u)}\left(\tau^{(\alpha)}_{B((\widetilde 0, R), R)}\leq t< \tau^{(\alpha)}_{H}\right)du\nn\\
&=&2^{d-1}|\partial D|\int_{0}^{R/2}
\Prob_{(\widetilde 0,u)}\left(t^{-1}\tau^{(\alpha)}_{B((\widetilde 0, R), R)}\leq 1< t^{-1}\tau^{(\alpha)}_{H}\right)du\nn\\
&=&2^{d-1}|\partial D|\int_{0}^{R/2}
\Prob_{(\widetilde 0, t^{-1/\alpha}u)}\left(\tau^{(\alpha)}_{B((\widetilde 0, t^{-1/\alpha}R), t^{-1/\alpha}R)}\leq 1< \tau^{(\alpha)}_{H}\right)du\nn\\
&=&2^{d-1}|\partial D|t^{1/\alpha}\int_{0}^{2^{-1}t^{-1/\alpha}R}
\Prob_{(\widetilde 0, v)}\left(\tau^{(\alpha)}_{B((\widetilde 0, t^{-1/\alpha}R), t^{-1/\alpha}R)}\leq 1< \tau^{(\alpha)}_{H}\right)dv\nn\\
&=&2^{d-1}|\partial D|t^{1/\alpha}\int_{0}^{\infty}1_{(0,2^{-1}t^{-1/\alpha}R)}(v)\Prob_{(\widetilde 0, v)}\left(\tau^{(\alpha)}_{B((\widetilde 0, t^{-1/\alpha}R), t^{-1/\alpha}R)}\leq 1< \tau^{(\alpha)}_{H}\right)dv\nn,
\end{eqnarray}
where we used the change of variables $v=t^{-1/\alpha}u$.

Now we will show that there exists a non-negative function $f$ on $(0, \infty)$ such that 
$$
1_{(0,R/2t^{1/\alpha})}(v)
\Prob_{(\widetilde 0, v)}\left(\tau^{(\alpha)}_{B((\widetilde 0, t^{-1/\alpha}R), 
t^{-1/\alpha}R)}\leq 1< \tau^{(\alpha)}_{H}\right)\leq f(v) \quad \text{and }\int_{0}^{\infty}f(v)dv<\infty.
$$

Assume $1\leq v$. 
Similarly as in the proof of Lemma \ref{lemma:inside}, by using \eqref{eqn:exp moment} and \eqref{eqn:BM exit} we find that 
\begin{eqnarray*}
&&\Prob_{(\widetilde 0, v)}\left(\tau^{(\alpha)}_{B((\widetilde 0, t^{-1/\alpha}R),
t^{-1/\alpha}R)}\leq 1< \tau^{(\alpha)}_{H}\right)
\leq \Prob_{(\widetilde 0, v)}\left( \tau_{B((\widetilde 0, v),v)}^{(\alpha)}\leq 1\right)
\leq c(\alpha,d)v^{-\alpha}.
\end{eqnarray*}
Hence, we have
$$
1_{(0,2^{-1}t^{-1/\alpha}R)}(v)
\Prob_{(\widetilde 0, v)}\left(\tau^{(\alpha)}_{B((\widetilde 0, t^{-1/\alpha}R), 
t^{-1/\alpha}R)}\leq 1< \tau^{(\alpha)}_{H}\right)\leq
1_{(0,1)}(v)+c(d,\alpha)1_{(1,\infty)}(v)\cdot v^{-\alpha}:=f(v).
$$
Starting from $(\widetilde 0, v)$ with $v\in (0,2^{-1}t^{-1/\alpha}R)$, as  $t\to 0$, 
$B((\widetilde 0, t^{-1/\alpha}R), t^{-1/\alpha}R)$ increases to $H$, and this implies $\tau^{(\alpha)}_{B((\widetilde 0, t^{-1/\alpha}R), t^{-1/\alpha}R)}\uparrow \tau^{(\alpha)}_{H}$.
Hence, we have
$$
\lim_{t\to 0}\Prob_{(\widetilde 0, v)}\left(\tau^{(\alpha)}_{B((\widetilde 0, t^{-1/\alpha}R), 
t^{-1/\alpha}R)}\leq 1< \tau^{(\alpha)}_{H}\right)=0.
$$
Hence, it follows from the Lebesgue dominated convergence theorem that
\begin{align*}
&\lim_{t\to 0}\frac{\int_{D\setminus D_{R/2}}\Prob_{x}\left(\tau_{D}^{(\alpha)}\leq t< \tau^{(\alpha)}_{H_{x}}\right)dx}{t^{1/\alpha}}\\
\leq&\lim_{t\to 0}|\partial D|\int_{0}^{\infty}1_{(0,2^{-1}t^{-1/\alpha}R)}(v)
\Prob_{(v,0)}\left(\tau^{(\alpha)}_{B((\widetilde 0, t^{-1/\alpha}R), 
t^{-1/\alpha}R)}\leq 1< \tau^{(\alpha)}_{H}\right)dv
=0.
\end{align*}
\qed

Now we establish the lower bound. 
\begin{lemma}\label{lemma:lb}
Let $D$ be a bounded $C^{1,1}$ open set in $\R^{d}$. Then, we have
$$
\liminf_{t\to 0}\frac{|D|-Q^{(\alpha)}_{D}(t)}{t^{1/\alpha}}\geq |\partial D|
\E[\overline{Y}^{(\alpha)}_{1}],
$$
where $\overline{Y}_{t}^{(\alpha)}=\sup\{Y_{u}^{(\alpha)}: 0\leq u\leq t\}$ and $Y^{(\alpha)}=\{Y_{t}^{(\alpha)}\}_{t\geq 0}$ is isotropic $\alpha$-stable process in $\R$.
\end{lemma}
\pf
Let $D$ satisfy the uniform interior and exterior $R$-ball condition.
Fix $a\leq R/2$. 
For $x\in D\setminus D_{a}$, let 
$z_x\in \partial D$ be such that $|x-z_x|=\delta_D(x)$ and  let ${\bf n}_{z_x}$ be the outward unit normal vector to the boundary $\partial D$ at the point $z_x$.

For $X$ starting from $x$, we define 
$Y^{(\alpha)}_{t}:=(X_{t}^{(\alpha)}-x)\cdot {\bf n}_{z_x}$, where $\cdot$ stands for usual scalar product in $\R^{d}$. Obviously, $Y^{(\alpha)}_{0}=0$.
Note that the characteristic exponent of $Y^{(\alpha)}_{t}$ is given by 
$$
\E[e^{i\eta Y^{(\alpha)}_{t}}]=
\E_x[e^{i\eta (X_{t}^{(\alpha)}-x)\cdot{\bf n}_{z_x}}]
=e^{-t|\eta|^{\alpha}}, 
\quad \eta \in \R,
$$
and this shows that $Y^{(\alpha)}_{t}$ is a one dimensional stable process starting from 0.
Let $r\le a$ and $x\in  \partial D_{r}$, and let $H_{x}$ be the half-space containing the interior $R$-ball at the point  $z_x$ and tangent to $\partial D$ at $z_x$.
When the process $X$ starts from $x$, we have
\begin{align*}
&\{\overline{Y}_{t}^{(\alpha)} >r\}=\{\tau_{H_x}^{(\alpha)}\leq t\} \subset \{\tau_{H_x}^{(\alpha)}\leq t, \tau_{D}^{(\alpha)} \leq t\} \cup \{\tau_{H_x}^{(\alpha)}\leq t<\tau_{D}^{(\alpha)} \}\\
\subset & \{\tau_{D}^{(\alpha)} \leq t\} \cup \{\tau_{H_x}^{(\alpha)}\leq t <\tau_{D}^{(\alpha)} \}
\subset  \{\tau_{D}^{(\alpha)} \leq t\} \cup \{\tau_{H_x}^{(\alpha)}\leq t <\tau_{\bar{B_x}^{c}}^{(\alpha)} \},
\end{align*}
where $\bar{B_x}$ is the unique exterior ball of radius $R$ in $D^{c}$ touching the point $z_x$ such that $\delta_{D}(x)=|x-z_{x}|$.
This implies 
$$
\Prob(\overline{Y}_{t}^{(\alpha)} >r)=\Prob_{x}(\tau_{H_x}^{(\alpha)}\leq t)\leq \Prob_{x}(\tau_{D}^{(\alpha)} \leq t) + \Prob_{x}(\tau_{H_x}^{(\alpha)}\leq t <\tau_{\bar{B_x}^{c}}^{(\alpha)} ).
$$
Hence, by the coarea formula
and \eqref{eqn:stability}, we have
\begin{align*}
&|D|-Q^{(\alpha)}_{D}(t)=\int_{D}\Prob_{x}(\tau^{(\alpha)}_{D}\leq t)dx
\geq \int_{D\setminus D_{a}}\Prob_{x}(\tau^{(\alpha)}_{D}\leq t)dx\\
\geq &\int_{D\setminus D_{a}}\Prob_{x}(\tau_{H_x}^{(\alpha)}\leq t)dx-\int_{D\setminus D_{a}}\Prob_{x}(\tau_{H_x}^{(\alpha)}\leq t <\tau_{\bar{B_x}^{c}}^{(\alpha)} )dx\\
\geq &\int_{D\setminus D_{a}}\Prob_{x}(\tau_{H_x}^{(\alpha)}\leq t)dx-\int_{D\setminus D_{R/2}}\Prob_{x}(\tau_{H_x}^{(\alpha)}\leq t <\tau_{\bar{B_x}^{c}}^{(\alpha)} )dx.
\end{align*}
It follows from Lemma \ref{lemma:inside} and Proposition \ref{prop:ub12} that
$$
\lim_{t\to0}t^{-1/\alpha}\int_{D\setminus D_{a}}\Prob_{x}(\tau_{H_x}^{(\alpha)}\leq t)dx = |\partial D|
\E[\overline{Y}^{(\alpha)}_{1}].
$$ 
Hence, it is enough to show that 
$$
\lim_{t\to 0}t^{-1/\alpha}\int_{D\setminus D_{R/2}}\Prob_{x}(\tau_{H_x}^{(\alpha)}\leq t <\tau_{\bar{B_x}^{c}}^{(\alpha)} )dx=0. 
$$

Now we use an idea similar to Lemma \ref{lemma:near boundary}. Note that by the rotational invariance 
of $X$ 
\begin{align}\label{eqn:approx1}
\int_{D\setminus D_{R/2}}\Prob_{x}\left(\tau_{H_x}^{(\alpha)}\leq t <\tau_{\bar{B_x}^{c}}^{(\alpha)} \right)dx
\leq \int_{0}^{R/2}|\partial D_{u}|\Prob_{(\widetilde 0, u)}\left( \tau^{(\alpha)}_{H} \leq t <\tau^{(\alpha)}_{B((\widetilde 0, -R), R)^{c}} \right)du,
\end{align}
where $H=\{x=(x_{1},\cdots ,x_{d}) : x_{d}>0\}$.
By the same argument as in Lemma \ref{lemma:near boundary}, we have 
\begin{align*}
&\int_{0}^{R/2}|\partial D_{u}|\Prob_{(\widetilde 0, u)}\left( \tau^{(\alpha)}_{H} \leq t <\tau^{(\alpha)}_{B((\widetilde 0, -R), R)^{c}} \right)du \\
\leq&
2^{d-1}|\partial D|t^{1/\alpha}\int_{0}^{\infty}1_{(0,R/2t^{1/\alpha})}(v)\Prob_{(\widetilde 0, v)}\left(\tau^{(\alpha)}_{H} \leq 1 <\tau^{(\alpha)}_{B((\widetilde 0, -R/t^{1/\alpha}), R/t^{1/\alpha})^{c}} \right)dv.
\end{align*}
We show that there exists a nonnegative function $g$ on $(0,\infty)$ such that
$$
1_{(0,R/2t^{1/\alpha})}(v)\Prob_{(\widetilde 0, v)}\left(\tau^{(\alpha)}_{H} \leq 1 <\tau^{(\alpha)}_{B((\widetilde 0, -R/t^{1/\alpha}), R/t^{1/\alpha})^{c}} \right)
\leq g(v) \text{ and } \int_{0}^{\infty}g(v)dv<\infty. 
$$
For $v\geq 1$, it follows from \eqref{eqn:exp moment} and \eqref{eqn:BM exit}
$$
\Prob_{(\widetilde 0, v)}\left(\tau^{(\alpha)}_{H} \leq 1 <\tau^{(\alpha)}_{B((\widetilde 0, -R/t^{1/\alpha}), R/t^{1/\alpha})^{c}} \right)
\leq
\Prob_{(\widetilde 0, v)}\left( \tau_{B((\widetilde 0, v),v)}^{(\alpha)}\leq 1\right)
\leq c(\alpha,d)v^{-\alpha}.
$$
Hence, we set
$$
g(v):=1_{(0,1)}(v)+c(d,\alpha)1_{(1,\infty)}(v)\cdot v^{-\alpha}.
$$
Since for each $v\in (0,R/2t^{1/\alpha})$ $\lim_{t\to 0}\Prob_{(\widetilde 0, v)}\left(\tau^{(\alpha)}_{H} \leq 1 <\tau^{(\alpha)}_{B((\widetilde 0, -R/t^{1/\alpha}), R/t^{1/\alpha})^{c}} \right)=0$, the conclusion follows from the Lebesgue dominated convergence theorem. 
\qed

\begin{theorem}\label{thm:main12}
Let $D\subset \R^{d}$, $d\geq 2$, be a bounded $C^{1,1}$ open 
set. 
Then, we have
\beq\label{eqn:main:12}
\lim_{t\to0}\frac{|D|-Q_{D}^{(\alpha)}(t)}{t^{1/\alpha}}=|\partial D|
\E[\overline{Y}^{(\alpha)}_{1}],
\eeq
where $\overline{Y}_{t}^{(\alpha)}$ is defined in \eqref{eqn:running}.
\end{theorem}
\pf
The lower bounded is proved in Lemma \ref{lemma:lb}.
Now we establish the upper bound. 
Assume $D$ satisfies the uniform interior and exterior $R$-ball condition.
For $\delta_{D}(x)<R/2$, it follows from \eqref{eqn:ub cases2} that
$$
\{\tau_{D}^{(\alpha)}\leq t\} \subset \{\tau_{H_{x}}^{(\alpha)}\leq t\} \cup \{\tau^{(\alpha)}_{D}\leq t<\tau^{(\alpha)}_{H_{x}}\},
$$
and 
$$
\int_{D\setminus D_{R/2}}\Prob_{x}(\tau_{D}^{(\alpha)}\leq t)\leq \int_{D\setminus D_{R/2}}\Prob_{x}(\tau_{H_{x}}^{(\alpha)}\leq t)dx + \int_{D\setminus D_{R/2}}\Prob_{x}(\tau^{(\alpha)}_{D}\leq t<\tau^{(\alpha)}_{H_{x}})dx.
$$
Combining inequality in \eqref{eqn:ub cases}, Lemmas  \ref{lemma:inside}, \ref{lemma:near boundary} and Proposition \ref{prop:ub12},  
we arrive at 
$$
\limsup_{t\to 0}\frac{|D|-Q_{D}^{(\alpha)}(t)}{t^{1/\alpha}}\leq |\partial D|
\E[\overline{Y}^{(\alpha)}_{1}],
$$
which together with Lemma \ref{lemma:lb} leads to the desired limit. 
\qed

\section{The case of Cauchy processes}\label{section:alpha1}
In this section, we establish 
the two-term asymptotic expansion for the spectral heat content of the Cauchy process.
Since $|D|-H_{D}^{(1)}(t)\leq |D|-Q_{D}^{(1)}(t)$, by \cite[Theorem 1.2 (ii)]{Val2017} we have
\beq\label{eqn:Cauchy lb}
\liminf_{t\to 0}\frac{|D|-Q_{D}^{(1)}(t)}{t\ln(1/t)}\geq\frac{|\partial D|}{\pi}.
\eeq
\begin{lemma}\label{lemma:exp sup}
Suppose $\alpha \in (1,2)$. It holds that
$$
\E[\overline{Y}^{(\alpha)}_{1}]
=\frac{\Gamma(1-\frac{1}{\alpha})}{\pi} +\int_{0}^{\infty}
\Prob(\overline{Y}_{1}^{(\alpha)}>u, Y_{1}^{(\alpha)}<u)du,
$$
where $\overline{Y}_{t}^{(\alpha)}$ is defined in \eqref{eqn:running}.
\end{lemma}
\pf
Consider the spectral heat content $Q_{(a,b)}^{(\alpha)}(t)$ of 
$(a,b)\subset \R$
for $Y_{t}^{(\alpha)}$. It holds that
\begin{eqnarray}\label{eqn:referee}
&&\frac{(b-a)-Q_{(a,b)}^{(\alpha)}(t)}{t^{1/\alpha}}
=\frac{\int_{a}^{b}\Prob_{x}(\tau_{(a,b)}^{(\alpha)}\leq t)dx}{t^{1/\alpha}}\nn\\
&=&\frac{\int_{a}^{b}\Prob_{x}(Y_{t}^{(\alpha)}\notin (a,b))dx}{t^{1/\alpha}}+\frac{\int_{a}^{b}\Prob_{x}(\tau_{(a,b)}^{(\alpha)}\leq t, Y_{t}^{(\alpha)}\in (a,b))dx}{t^{1/\alpha}}.
\end{eqnarray}
We will show that 
\beq\label{eqn:Cauchy aux1}
\lim_{t\to 0}\frac{\int_{a}^{b}\Prob_{x}(\tau_{(a,b)}^{(\alpha)}\leq t, 
Y_{t}^{(\alpha)}\in (a,b))dx}{t^{1/\alpha}}=2\int_{0}^{\infty}\Prob(\overline{Y}_{1}^{(\alpha)}>u, Y_{1}^{(\alpha)}<u)du.
\eeq
Together with $\lim_{t\to 0}t^{-1/\alpha}\int_{a}^{b}\Prob_{x}(Y_{t}^{(\alpha)}\notin (a,b))dx=2\pi^{-1}\Gamma(1-\frac{1}{\alpha})$ from \cite[Theorem 1.1. (i)]{Val2017} and
$\lim_{t\to 0}t^{-1/\alpha}\int_{a}^{b}\Prob_{x}(\tau_{(a,b)}^{(\alpha)}\leq t)dx=2\E[\overline{Y}^{(\alpha)}_{1}]$ from \cite[Theorem 1.1]{Val2016}, 
this will establish the conclusion by taking limits at both sides of \eqref{eqn:referee}. 

Note that we have 
\begin{eqnarray*}
&&\int_{a}^{b}\Prob_{x}(\tau_{(a,b)}^{(\alpha)}\leq t, Y_{t}^{(\alpha)}\in (a,b))dx\nn\\
&=&\int_{a}^{b}\Prob_{x}(\overline{Y}^{(\alpha)}_{t} >b, Y_{t}^{(\alpha)}\in (a,b))dx +\int_{a}^{b}\Prob_{x}(\underline{Y}^{(\alpha)}_{t}<a, Y_{t}^{(\alpha)}\in (a,b))dx\nn\\
&&-\int_{a}^{b}\Prob_{x}(\overline{Y}^{(\alpha)}_{t} >b, \underline{Y}^{(\alpha)}_{t}<a, Y_{t}^{(\alpha)}\in (a,b))dx,
\end{eqnarray*}
where $\underline{Y}^{(\alpha)}_{t}=\inf\{ Y_{u}^{(\alpha)}: 0\leq u\leq t\}$ is the running infimum process of $Y_{t}^{(\alpha)}$.
It follows from \cite[Lemma 3.2]{P20} that
$$
\int_{a}^{b}\Prob_{x}(\overline{Y}^{(\alpha)}_{t} >b, \underline{Y}^{(\alpha)}_{t}<a, Y_{t}^{(\alpha)}\in (a,b))dx=O(t^{1+\frac{1}{\alpha}}) \text{ as } t\downarrow 0.
$$
By the change of variables $u=(b-x)t^{-1/\alpha}$ and the scaling property we have
\begin{eqnarray*}
&&\int_{a}^{b}\Prob_{x}(\overline{Y}^{(\alpha)}_{t} >b, Y_{t}^{(\alpha)}\in (a,b))dx\\
&=&t^{1/\alpha}\int_{0}^{(b-a)t^{-1/\alpha}}\Prob(\overline{Y}^{(\alpha)}_{1}>u, -(b-a)t^{-1/\alpha}+u <Y_{1}^{(\alpha)}<u)du.
\end{eqnarray*}
Note that we have
$$
\Prob(\overline{Y}^{(\alpha)}_{1}>u, -(b-a)t^{-1/\alpha}+u <Y_{1}^{(\alpha)}<u)\leq \Prob(\overline{Y}^{(\alpha)}>u)
$$
and
$$
\int_{0}^{\infty}\Prob(\overline{Y}^{(\alpha)}_{1}>u)du=\E[\overline{Y}^{(\alpha)}_{1}]<\infty.
$$
Hence, it follows from the Lebesgue dominated convergence theorem that the limit is
\begin{align*}
&\lim_{t\to 0}\frac{\int_{a}^{b}\Prob_{x}(\tau_{(a,b)}^{(\alpha)}\leq t, Y_{t}^{(\alpha)}\in (a,b))dx}{t^{1/\alpha}}\\
=&\lim_{t\to 0}\frac{t^{1/\alpha}\int_{0}^{(b-a)t^{-1/\alpha}}\Prob(\overline{Y}^{(\alpha)}_{1}>u, -(b-a)t^{-1/\alpha}+u <Y_{1}^{(\alpha)}<u)du}{t^{1/\alpha}}\\
=&\int_{0}^{\infty}\Prob(\overline{Y}^{(\alpha)}_{1}>u, Y_{1}^{(\alpha)}<u)du.
\end{align*}
The second term can be handled in a similar way using the symmetry of $Y^{(\alpha)}_{t}$, and this proves \eqref{eqn:Cauchy aux1}.
\qed

We need a simple lemma which is similar to \cite[Lemma 3.2]{PS19}. The proof is essentially the same with obvious modifications. We provide the details for the reader's convenience. 
\begin{lemma}\label{lemma:SS expectation}
For any $\delta>0$ and $\alpha\in (0,2)$, we have
$$
\lim_{t\to 0}\frac{\E[\left(S_{1}^{(\alpha/2)}\right)^{\alpha/2}, 0<S_{1}^{(\alpha/2)}<\delta t^{-2/\alpha} ]}{\ln(1/t)}=\frac{1}{\Gamma(1-\frac{\alpha}{2})}.
$$
\end{lemma}
\pf
Note that the stable subordinator $S_{t}^{(\alpha/2)}$ has a continuous transition density $g^{(\alpha/2)}(t,u)$ and it follows from \cite[(2.5)]{PS19} that
\beq\label{eqn:SS limit}
\lim_{u\to\infty}g^{(\alpha/2)}(1,u)u^{1+\frac{\alpha}{2}}=\frac{\alpha}{2\Gamma(1-\frac{\alpha}{2})}.
\eeq
Also it follows from \cite[Proposition 2.1]{Val2017} that 
$$
\lim_{t\to 0}\E[\left(S_{1}^{(\alpha/2)}\right)^{\alpha/2}, 0<S_{1}^{(\alpha/2)}<\delta t^{-2/\alpha} ]=\infty.
$$
Hence, it follows from L'H\^opital's rule and \eqref{eqn:SS limit} that
\begin{eqnarray*}
&&\lim_{t\to 0}\frac{\E[\left(S_{1}^{(\alpha/2)}\right)^{\alpha/2}, 0<S_{1}^{(\alpha/2)}<\delta t^{-2/\alpha} ]}{\ln(1/t)}=\lim_{t\to 0}\frac{\int_{0}^{\delta t^{-2/\alpha}} u^{\alpha/2}g^{(\alpha/2)}(1,u)du}{\ln(1/t)}\\
&=&\lim_{t\to 0}\frac{(\delta t^{-2/\alpha})^{\alpha/2}g^{(\alpha/2)}(1,\delta t^{-2/\alpha})(-\frac{2}{\alpha}\delta t^{-\frac{2}{\alpha}-1}) }{-1/t}=\lim_{t\to0}\frac{2}{\alpha}g^{(\alpha/2)}(1,\delta t^{-2/\alpha})(\delta t^{-2/\alpha})^{1+\frac{\alpha}{2}}\\
&=&\frac{1}{\Gamma(1-\frac{\alpha}{2})}.
\end{eqnarray*}
\qed

\begin{lemma}\label{lemma:Cauchy ub}
Let $D$ be a bounded $C^{1,1}$ open set in $\R^{d}$, $d\geq 2$. For any $\alpha\in (1,2)$, we have
$$
\limsup_{t\to 0}\frac{|D|-Q_{D}^{(1)}(t)}{t\ln(1/t)}\leq \frac{|\partial D|}{\pi} +
\frac{\int_{0}^{\infty}\Prob(
\overline{Y}^{(\alpha)}_{1}>u, Y_{1}^{(\alpha)}<u )du }
{\Gamma(1-\frac{1}{\alpha})}\cdot |\partial D|,
$$
where $\overline{Y}_{t}^{(\alpha)}$ is defined in \eqref{eqn:running}.
\end{lemma}
\pf
The main tool for this proof is the spectral heat content for \textit{subordinate killed stable processes}. 
Let $\alpha\in(1,2)$ and $\alpha\beta=2$. Let $S_{t}$ and $T_{t}$ be independent $\frac{\alpha}{2}$ and $\frac{\beta}{2}$ stable subordinators which are independent of the Brownian motion $W_{t}$. 
Let $X^{(\alpha)}_t=W_{S_{t}^{(\alpha/2)}}$, $X^{(\alpha), D}_t$ be the process $X^{(\alpha)}_t$ killed upon exiting $D$, and
$Z_{t}=X^{(\alpha), D}_{T_{t}^{(\beta/2)}}$ be 
the subordinate killed $\alpha$-stable process.
Let $\widetilde{Q}^{(\alpha,\beta)}_{D}(t)$ be the spectral heat content for $Z$, that is, 
$$
\widetilde{Q}^{(\alpha,\beta)}_{D}(t)=\int_{D}\Prob_{x}(\tau_{D}^{(\alpha)}>T_{t}^{(\beta/2)})dx.
$$

Since $\{\tau_{D}^{(1)} \leq t\}\subset \{\tau_{D}^{(\alpha)} \leq T_{t}^{(\beta/2)}\}$, we have $|D|-Q_{D}^{(1)}(t)\leq |D|-\widetilde{Q}^{(\alpha,\beta)}_{D}(t)$. 
We will show that 
\beq\label{eqn:Cauchy aux7}
\limsup_{t\to 0}\frac{|D|-\widetilde{Q}^{(\alpha,\beta)}_{D}(t)}{t\ln(1/t)}\leq 
\frac{|\partial D|}{\pi} +
\frac{\int_{0}^{\infty}
\Prob(\overline{Y}^{(\alpha)}_{1}>u, Y_{1}^{(\alpha)}<u )du }
{\Gamma(1-\frac{1}{\alpha})}\cdot |\partial D|.
\eeq
Note that 
\begin{eqnarray}\label{eqn:Cauchy aux3}
&&|D|-\tilde{Q}^{(\alpha,\beta)}_{D}(t)=\int_{0}^{\infty}(|D|-Q_{D}^{(\alpha)}(u))\Prob(T_{t}^{(\beta/2)}\in du).
\end{eqnarray}
It follows from Lemma \ref{lemma:exp sup} and \eqref{eqn:main:12} that, for any $\eps>0$,
there exists $\delta=\delta(\eps)$ such that for all $u\leq \delta$,
\begin{eqnarray*}
\frac{|D|-Q_{D}^{(\alpha)}(u)}{u^{1/\alpha}} 
\leq \left(\frac{\Gamma(1-\frac{1}{\alpha})}{\pi} + \int_{0}^{\infty}\Prob(\overline{Y}^{(\alpha)}_{1}>u, Y_{1}^{(\alpha)}<u )du\right)|\partial D| +\eps.
:=C_{1}+\eps.
\end{eqnarray*}
Hence, \eqref{eqn:Cauchy aux3} can be written as
\begin{eqnarray}\label{eqn:Cauchy aux4}
&&\int_{0}^{\infty}(|D|-Q_{D}^{(\alpha)}(u))\Prob(T_{t}^{(\beta/2)}\in du)\nn\\
&=&\int_{0}^{\delta}(|D|-Q_{D}^{(\alpha)}(u))\Prob(T_{t}^{(\beta/2)}\in du)+\int_{\delta}^{\infty}(|D|-Q_{D}^{(\alpha)}(u))\Prob(T_{t}^{(\beta/2)}\in du)\nn\\
&=&\int_{0}^{\delta}(\frac{|D|-Q_{D}^{(\alpha)}(u)}{u^{1/\alpha}})u^{1/\alpha}\Prob(T_{t}^{(\beta/2)}\in du)+\int_{\delta}^{\infty}(|D|-Q_{D}^{(\alpha)}(u))\Prob(T_{t}^{(\beta/2)}\in du)\nn\\
&\leq&\int_{0}^{\delta}( C_{1}+\eps)u^{1/\alpha}\Prob(T_{t}^{(\beta/2)}\in du)+\int_{\delta}^{\infty}(|D|-Q_{D}^{(\alpha)}(u))\Prob(T_{t}^{(\beta/2)}\in du).
\end{eqnarray}
It follows from \cite[(2.8)]{PS19} that the second term in 
\eqref{eqn:Cauchy aux4} can be estimated by
\begin{eqnarray*}
&&\int_{\delta}^{\infty}(|D|-Q_{D}^{(\alpha)}(u))\Prob(T_{t}^{(\beta/2)}\in du)\leq|D|\int_{\delta}^{\infty}\Prob(t^{2/\beta}T_{1}^{(\beta/2)}\in du)\\
&=&|D|\int_{\delta t^{-2/\beta}}^{\infty}\Prob(T_{1}^{(\beta/2)}\in dv)\leq c|D|(\delta t^{-2/\beta})^{-\frac{\beta}{2}}=c|D|\delta^{-\frac{\beta}{2}}t,
\end{eqnarray*}
and 
\beq\label{eqn:Cauchy aux5}
\limsup_{t\to 0}\frac{\int_{\delta}^{\infty}(|D|-Q_{D}^{(\alpha)}(u))\Prob(T_{t}^{(\beta/2)}\in du)}{t\ln(1/t)}\leq \limsup_{t\to 0}\frac{c|D|\delta^{-\frac{\beta}{2}}t}{t\ln(1/t)}=0.
\eeq

Recall that $\alpha\beta=2$. 
Hence, by the change of variables $u=t^{2/\beta}v$, the first expression in \eqref{eqn:Cauchy aux4} can be written as
\begin{eqnarray*}
&&\int_{0}^{\delta}u^{1/\alpha}\Prob(T_{t}^{(\beta/2)}\in du)\nn=\int_{0}^{\delta}u^{1/\alpha}\Prob(t^{2/\beta}T_{1}^{(\beta/2)}\in du)\nn\\
&=&\int_{0}^{\delta t^{-2/\beta}}(t^{2/\beta}v)^{1/\alpha}\Prob(T_{1}^{(\beta/2)}\in dv)\nn=t\int_{0}^{\delta t^{-2/\beta}}v^{\beta/2}\Prob(T_{1}^{(\beta/2)}\in dv).
\end{eqnarray*}
It follows from Lemma \ref{lemma:SS expectation} that
\beq\label{eqn:Cauchy aux6}
\limsup_{t\to 0}\frac{t\int_{0}^{\delta t^{-2/\beta}}v^{\beta/2}\Prob(S_{1}^{(\beta/2)}\in dv)}{t\ln(1/t)}=\frac{1}{\Gamma(1-\frac{\beta}{2})}.
\eeq

Combining \eqref{eqn:Cauchy aux5} and \eqref{eqn:Cauchy aux6} we conclude that 
$$
\limsup_{t\to 0}\frac{|D|-\widetilde{Q}^{\alpha,\beta}(t)}{t\ln(1/t)}\leq ( C_{1} +\eps) \times \frac{1}{\Gamma(1-\frac{\beta}{2})}=\frac{|\partial D|}{\pi} 
+\frac{\int_{0}^{\infty}
\Prob(\overline{Y}^{(\alpha)}_{1}>u, Y_{1}^{(\alpha)}<u )du \cdot |\partial D|}
{\Gamma(1-\frac{1}{\alpha})}
+\frac{\eps}{\Gamma(1-\frac{1}{\alpha})}.
$$
Since $\eps>0$ is arbitrary, this establishes \eqref{eqn:Cauchy aux7}.
\qed

We now recall the definition for the double gamma function from \cite[(4.4)]{Ku11}. 
For $z\in \C$ and $\tau \in \C$ with $|\text{arg}(\tau)|<\pi$, we define 
$$
G(z;\tau)=\frac{z}{\tau}e^{\frac{az}{\tau} + \frac{bz^2}{2\tau}}
{\textstyle \prod}_{m\geq 0}{\textstyle \prod}^{'}_{n\geq 0}\left( 1+\frac{z}{m\tau +n}\right)\exp\left(-\frac{z}{m\tau +n} +\frac{z^2}{2(m\tau +n)^2}\right),
$$
where the prime in the second product means that the term corresponding to $m=n=0$ is omitted. 

\begin{lemma}\label{lemma:double gamma bounded}
Let $K$ be a compact set in $\C$. 
Then there exists a constant $c=c(K)>0$ such that 
$$
\left| G(z;\tau) \right| \leq c \text{ for all } z\in K \text{ and } \tau\in [1,2].
$$
\end{lemma}
\pf
On the compact set $K\times [1, 2]$, $\frac{z}{\tau}e^{\frac{az}{\tau} + \frac{bz^2}{2\tau}}$ is continuous and hence bounded. 
We only need to prove that the double infinite product in $G(z, \tau)$ is bounded. 
Recall the canonical factor $E_{2}(z)$ (see \cite[p. 145]{SS complex}):
$$
E_{2}(z)=(1-z)e^{z+\frac{z^2}{2}}.
$$
The double infinite product in $G(z;\tau)$ can be written as
$$
 {\textstyle \prod}_{m\geq 0}{\textstyle \prod}^{'}_{n\geq 0}\left( 1+\frac{z}{m\tau +n}\right)\exp\left(-\frac{z}{m\tau +n} +\frac{z^2}{2(m\tau +n)^2}\right)
 = {\textstyle \prod}_{m\geq 0}{\textstyle \prod}^{'}_{n\geq 0}E_{2}(-\frac{z}{m\tau +n}).
$$
It follows from \cite[Lemma 4.2]{SS complex} for $|w|\leq 1/2$, $|1-E_{2}(w)|\leq c_{1}|w|^{3}$ for some constant $c_{1}$.
Since $|-\frac{z}{m\tau +n}| \leq \frac{|z|}{m+n}$, 
all but finitely many terms of the form $|-\frac{z}{m\tau +n}|$ in the double infinite product are less than $\frac12$, and since each of these finitely many terms are bounded on $K\times [1, 2]$, we may assume $|-\frac{z}{m\tau +n}|\leq \frac12$ on $K\times [1, 2]$ for all $n,m$.

Now it follows from \cite[Lemma 4.2]{SS complex} for all 
$(z, \tau)\in  K\times [1, 2]$ that
$$
\left| 1-E_{2}(-\frac{z}{m\tau +n})\right| \leq c_{1}\left |\frac{z}{m\tau +n}\right|^{3}\leq \frac{c_{2}}{(m+n)^{3}} \text{ for all } n,m.
$$
It is easy to see that 
$$
\sum_{m\geq 0, n\geq 0, (m,n)\neq (0, 0)}\frac{1}{(m+n)^{3}}<\infty.
$$
Note that for any $|w|\leq \frac{1}{2}$ we have $\ln |1-w|\leq 2|w|$.
There exist positive integers $N_1$ and  $M_1$ such that $|1-E_{2}(-\frac{z}{m\tau +n})|\leq \frac{1}{2}$ on $K\times [1, 2]$ for all $n\geq N_{1}$ and $m\geq M_{1}$.
Hence, for any $M, N \in \N$ large we have 
\begin{align*}
&\left|{\textstyle \prod}_{M_{1}\leq m \leq M}{\textstyle \prod}_{N_{1}\leq n\leq N} E_{2}(-\frac{z}{m\tau +n}) \right|
={\textstyle \prod}_{M_{1}\leq m \leq M}{\textstyle \prod}_{N_{1}\leq n\leq N}\left|1- (1- E_{2}(-\frac{z}{m\tau +n}) )\right|\\
&={\textstyle \prod}_{M_{1}\leq m \leq M}{\textstyle \prod}_{N_{1}\leq n\leq N}\exp \left( \ln  \left| (1- (1- E_{2}(-\frac{z}{m\tau +n}) )) \right| \right)\\
&= \exp \left( \sum_{M_{1}\leq m \leq M,N_{1}\leq n\leq N} \ln \left|(1- (1- E_{2}(-\frac{z}{m\tau +n}) )) \right| \right)\\
&\leq \exp \left( 2\sum_{M_{1}\leq m \leq M,N_{1}\leq n\leq N} \left| 1- E_{2}(-\frac{z}{m\tau +n})  \right| \right)\leq \exp \left( 2c_2\sum_{M_{1}\leq m \leq M,N_{1}\leq n\leq N} \frac{1}{(m+n)^{3}} \right)\\
&\leq \exp \left( 2c_2\sum_{M_{1}\leq m, N_{1}\leq n} \frac{1}{(m+n)^{3}} \right).
\end{align*}
By letting $M, N\to \infty$ we see that the double infinite product is bounded on 
$K\times [1, 2]$.
\qed

It follows from \cite[Proposition VIII.1-4]{Ber} that there exists a constant $C>0$ such that
\beq\label{eqn:same asymptotic}
\Prob(\overline{Y}^{(\alpha)}_1>u)\sim\Prob(Y^{(\alpha)}_1>u)\sim Cu^{-\alpha} 
\text{ as } u\to \infty.
\eeq
Let $\rho:=\Prob(Y_1^{(\alpha)}>0)$.  We say $Y^{(\alpha)}\in C_{k,l}$ if
\begin{align}\label{eqn:ckl}
&(\alpha, \rho)\in \{\alpha\in (0,1), \rho\in(0,1)\} \cup \{\alpha=1, \rho=1/2\}\cup \{\alpha\in (1,2), \rho\in [1-\frac{1}{\alpha}, \frac{1}{\alpha}]\} \text{ and}\nn\\
&\frac{1}{2}+k=\frac{l}{\alpha}, \text{ or equivalently } \alpha=\frac{2l}{1+2k} \text{ for some } k,l\in \N.
\end{align}
Note that this condition already appeared in \cite[\text{Definition 1}]{Ku11}. 

\begin{lemma}\label{lemma:uniform ub}
Suppose that $Y^{(\alpha)}\in C_{k,k+1}$
for some $k\in \mathbb{N}$.
Then $\overline{Y}_1^{(\alpha)}$ has a density $\overline{p}^{(\alpha)}(x)$ and there exists 
a constant $A>0$, 
independent of all $\alpha\in(1,2)$ satisfying $Y_{t}^{(\alpha)}\in C_{k,k+1}$, such that
$$
\overline{p}^{(\alpha)}(x)\leq C\alpha x^{-1-\alpha} +Ax^{-3}
$$
for all $x>0$, where $C$ is the constant in \eqref{eqn:same asymptotic}.
\end{lemma}
\pf
The proof is similar to that of \cite[Theorem 9]{Ku11} with a focus on establishing 
a uniform constant $A$. 
It follows from the proof of \cite[Theorem 9]{Ku11} that $\overline{Y}_1^{(\alpha)}$ has a density $\overline{p}^{(\alpha)}(x)$ given by the inverse Mellin transform
$$
\overline{p}^{(\alpha)}(x)=\frac{1}{2\pi i}\int_{1+i\R}M(s,\alpha)x^{-s}ds,
$$
where $M(s,\alpha)$ is the Mellin transform given by
$$
M(s,\alpha)=\E[(\overline{Y}^{(\alpha)}_{1})^{s-1}], \quad s\in \mathbb{C}.
$$
We remark here that we write the Mellin transform as $M(s,\alpha)$ instead of $M(s)$ to emphasize its dependence on $\alpha$.
It follows from \cite[Lemma 2]{Ku11} that $M(s,\alpha)$ can be extended to a meromorphic function on $\mathbb{C}$ whose simples poles are at $s_{m,n}:=m+\alpha n$, where $m\leq 1-(k+1)$ and $n\in \{0,1,\cdots, k\}$, or $m\geq 1$ and $n\in \{1,2,\cdots, k\}$ with residues
$$
\text{Res}(M(s,\alpha), s_{m,n})=c_{m-1,n}^{+},
$$
where $c^{+}_{m-1,n}$ is the constant defined in \cite[(7.5)]{Ku11}.

Since $Y^{(\alpha)}\in C_{k,k+1}$, 
it follows from \eqref{eqn:ckl} that $\alpha=\frac{2+2k}{1+2k}\in (1,2)$.
If $m\geq 1$ and $n\in \{1,\cdots, k\}$ then $s_{m,n}=m+\alpha n$ has a nonempty intersection with  $[1,3]$ if and only if $m=n=1$. In particular, $s_{1,1}=1+\alpha$. 
If $m\leq 1-(k+1)=-k$ and $n\in \{0,\cdots, k\}$, we have 
$$
s_{m,n}=m+\alpha n \leq -k +\frac{2+2k}{1+2k}k=\frac{k}{1+2k} <\frac{1}{2}.
$$
Hence, the only simple pole of $M(s,\alpha)$ with $\text{Re}(s)\in [1,3]$ is $s_{1,1}=1+\alpha$.

We note here that although it is written as $y\to \infty$ in \cite[Lemma 3]{Ku11}, it is actually true for all $|y|\to \infty$ by checking the proof there as \cite[(7.15)]{Ku11} is true for all $z\to\infty$ with $|\text{arg}(z)|<\pi$.
Hence, by taking a rectangle with vertices at $1\pm Pi$ and $3\pm Pi$ with the help  of the residue theorem, then letting $P\to \infty$ and using \cite[Lemma 3]{Ku11},  we get
\begin{eqnarray*}
&&\overline{p}^{(\alpha)}(x)=\frac{1}{2\pi i}\int_{1+i\R}M(s,\alpha)x^{-s}ds \\
&=&\text{Res}(M(s,\alpha), s_{1,1}^{+})x^{-1-\alpha} +\frac{1}{2\pi i}\int_{3+i\R}M(s,\alpha)x^{-s}ds\\
&=&c_{0,1}^{+}x^{-1-\alpha}+\frac{1}{2\pi i}\int_{3+i\R}M(s,\alpha)x^{-s}ds,
\end{eqnarray*}
where $\text{Res}(M(s,\alpha), s_{1,1}^{+})$ is the residue of $M(s,\alpha)$ at $s_{1,1}^{+}$.

We claim that the constant $c_{0,1}^{+}$ must be $C\alpha$, where $C$ is from \eqref{eqn:same asymptotic} as we will show that the reminder is $O(x^{-3})$, which will in turn imply that there exist constants $c_{1}, c_{2}\in \R$ such that 
$$
c_{0,1}^{+}x^{-1-\alpha} +c_1x^{-3} \leq \overline{p}^{(\alpha)}(x) \leq c_{0,1}^{+}x^{-1-\alpha} +c_2x^{-3}
$$
for all sufficiently large $x$.
By integrating on $(u,\infty)$ we obtain 
$$
\frac{c_{0,1}^{+}}{\alpha}u^{-\alpha} +\frac{c_1}{2} u^{-2}\leq 
\Prob(\overline{Y}^{(\alpha)}_{1}>u) 
\leq \frac{c_{0,1}^{+}}{\alpha}u^{-\alpha} +\frac{c_2}{2}u^{-2},
$$
for all sufficiently large $u>0$.
Comparing the equation above with \eqref{eqn:same asymptotic} we conclude that 
$\alpha^{-1}c_{0,1}^{+}=C$.
This shows that the leading term of $\overline{p}^{(\alpha)}(x)$ is $C\alpha x^{-1-\alpha}$.

Hence, the proof will be complete once we show that 
\beq\label{eqn:Mellin error}
\left| \frac{1}{2\pi i}\int_{3+i\R}M(s,\alpha)x^{-s}ds\right| 
\leq Ax^{-3},
\eeq
where the constant $A$ is independent of all $\alpha\in (1,2)$ satisfying 
$Y^{(\alpha)}\in C_{k,k+1}$.

From \cite[Lemma 3]{Ku11} we know that 
\beq\label{eqn:Mellin}
\ln|M(3+iy,\alpha)|=-\frac{\pi |y|}{2\alpha} +o(y) \text{ as } |y|\to\infty.
\eeq
The proof of \eqref{eqn:Mellin} in \cite[Lemma 3]{Ku11} depends on \cite[(4.5)]{BK97}, which is a uniform estimate of the double gamma functions. Hence, the error term in \eqref{eqn:Mellin} is uniform for all $\alpha\in (1,2)$ and 
there exists $N>0$, independent of  $\alpha\in (1,2)$, such that 
$$
\ln|M(3+iy,\alpha)|\leq -\frac{\pi |y|}{5} \text{ for all } |y|\geq N \text{ and } \alpha\in(1,2).
$$
Hence, we have 
\beq\label{eqn:Mellin ub}
|M(3+iy, \alpha)|\leq e^{-\frac{\pi|y|}{5}} \text{ for all } |y|\geq N \text{ and } \alpha\in(1,2).
\eeq
It follows from \cite[Theorem 8]{Ku11} that the Mellin transform $M(s,\alpha)$ can be written as
$$
M(s,\alpha)=\alpha^{s-1}\frac{G(\alpha/2;\alpha)}{G(\alpha/2+1;\alpha)}\frac{G(\alpha/2+2-s;\alpha)}{G(\alpha/2-1+s;\alpha)}\frac{G(\alpha-1+s;\alpha)}{G(\alpha+1-s;\alpha)},
$$
where $G(z;\tau)$ is the double gamma function.
$G(z;\tau)$ has simple zeroes on the lattice $m\tau+n$, $m,n\leq 0$.
By a simple calculation one can check that the double gamma functions in the denominators above
have no zero with $s=3+iy$ for $|y|\leq N$.
It follows from Lemma \ref{lemma:double gamma bounded} that there exists a constant $c>0$ such that
$$
|M(3+iy,\alpha)|\leq ce^{-\frac{\pi|y|}{5}} \text{ for all } |y|\leq N \text{ and } \alpha \in (1,2),
$$
and combining with \eqref{eqn:Mellin ub} we have 
$$
|M(3+iy,\alpha)|\leq ce^{-\frac{\pi|y|}{5}} \text{ for all } y\in \R \text{ and } \alpha \in (1,2).
$$

Hence, \eqref{eqn:Mellin error} can be estimated as 
$$
\left| \frac{1}{2\pi i}\int_{3+i\R}M(s,\alpha)x^{-s}ds\right| \leq \frac{x^{-3}}{2\pi}\int_{-\infty}^{\infty}|M(3+it, \alpha)|dt\leq\frac{x^{-3}}{2\pi}\int_{-\infty}^{\infty}ce^{-\frac{\pi |t|}{5}}dt:=
Ax^{-3}.
$$
\qed

Let $p^{(\alpha)}(x)$ be the transition density of $Y^{(\alpha)}_1$.
The following lemma is a variation of \cite[Propositions 7.1.1 and 7.1.2]{Ko} with an explicit error estimate. 
\begin{lemma}\label{lemma:uniform lb}
Let $\alpha\in(1,2)$ and $n\in\N$. Then, we have 
$$
p^{(\alpha)}(x)=
\frac{1}{\pi}\sum_{k=1}^{n}\frac{(-1)^{k+1}\Gamma(1+k\alpha)\sin(\frac{k\alpha\pi}{2}) }{k!}x^{-1-\alpha k}+E(x,\alpha),
$$
where 
$$
\left|E(x,\alpha)\right|\leq \frac{2\Gamma(2n+2)}{\pi n!}x^{-2-n} \text{ for all } x\geq 1 \text{ and } \alpha \in(1,2).
$$
In particular, for any $x\geq 1$ and $\alpha\in (1,2)$ we have
$$
 \frac{\Gamma(1+\alpha)\sin\frac{\pi\alpha}{2}}{\pi}x^{-1-\alpha}-\frac{12}{\pi}x^{-3} \leq p^{(\alpha)}(x)\leq \frac{\Gamma(1+\alpha)\sin\frac{\pi\alpha}{2}}{\pi}x^{-1-\alpha}+\frac{12}{\pi}x^{-3},
$$
and in fact 
$$
 \frac{\Gamma(1+\alpha)\sin\frac{\pi\alpha}{2}}{\pi}=C\alpha,
$$
where $C$ is the constant in \eqref{eqn:same asymptotic}.
\end{lemma}
\pf
To prove this lemma, we consider a 1-dimensional $\alpha$-stable processes $Y^{(\alpha), \gamma}_{t}$ with
$$
\E[e^{iy Y_{t}^{(\alpha), \gamma}}]=\exp(-t|y|^{\alpha}e^{\frac{i\pi\gamma}{2}\text{sgn}(y)}),
$$
where $\gamma$ represents the skewness of $Y^{(\alpha), \gamma}_{t}$  and $\text{sgn}(y):=1_{\{y\geq 0\}}- 1_{\{y<0\}}$. 
We use 
$p^{(\alpha)}(x, \gamma)$
to denote the density of $Y^{(\alpha), \gamma}_{1}$. When $\gamma=0$, $Y^{(\alpha), \gamma}_{t}$ reduces to $Y^{(\alpha)}_{t}$.
By the Fourier inversion we have 
\begin{eqnarray*}
&&p^{(\alpha)}(x,\gamma)
=\frac{1}{2\pi}\int_{-\infty}^{\infty}e^{-ixy}\exp(-|y|^{\alpha}e^{\frac{i\pi\gamma}{2}\text{sgn}(y)})dy\\
&=&\frac{1}{2\pi}\int_{0}^{\infty}e^{-ixy}\exp(-|y|^{\alpha}e^{\frac{i\pi\gamma}{2}})dy+\frac{1}{2\pi}\int_{-\infty}^{0}e^{-ixy}\exp(-|y|^{\alpha}e^{-\frac{i\pi\gamma}{2}})dy\\
&=&\frac{1}{2\pi}\int_{0}^{\infty}e^{-ixy}\exp(-|y|^{\alpha}e^{\frac{i\pi\gamma}{2}})dy+\frac{1}{2\pi}\int^{\infty}_{0}e^{ixy}\exp(-|y|^{\alpha}e^{-\frac{i\pi\gamma}{2}})dy\\
&=&\frac{1}{\pi}\text{Re}\int_{0}^{\infty}e^{-ixy}\exp(-y^{\alpha}e^{\frac{i\pi\gamma}{2}})dy,
\end{eqnarray*} 
where $\text{Re}(\cdot)$ represents the real part of the argument. 
By the Taylor expansion with remainder we can write $e^{-ixy}$ as
$$
e^{-ixy}=\sum_{k=0}^{n-1}\frac{(-ix)^{k}}{k!}y^{k} +\frac{x^{n}y^{n}}{n!}R,
$$
where $|R|\leq 1$.
Since, 
$$
\int_{0}^{\infty}y^{\beta}\exp(-\lam y^{\alpha})dy=\alpha^{-1}\lam^{-\frac{\beta+1}{\alpha}}\Gamma(\frac{\beta+1}{\alpha}),
$$
we have 
\begin{align}\label{eqn:HK}
&p^{(\alpha)}(x,\gamma)
=\frac{1}{\pi}\text{Re}\int_{0}^{\infty}\left( \sum_{k=0}^{n-1}\frac{(-ix)^{k}}{k!}y^{k} +\frac{x^{n}y^{n}}{n!}R\right)\exp(-y^{\alpha}e^{\frac{i\pi\gamma}{2}})dy\nn\\
&=\frac{1}{\pi}\text{Re}\sum_{k=0}^{n-1}\frac{(-ix)^{k}}{k!}\frac{1}{\alpha}(e^{\frac{i\gamma\pi}{2}})^{-\frac{k+1}{\alpha}}\Gamma(\frac{k+1}{\alpha})+
\text{Re}\left(\frac{x^{n}}{\pi\alpha n!}\Gamma(\frac{n+1}{\alpha})(e^{\frac{i\gamma\pi}{2}})^{-\frac{n+1}{\alpha}}R\right)\nn\\
&=\frac{1}{\alpha\pi}\text{Re}\sum_{k=1}^{n}\frac{(-ix)^{k-1}}{(k-1)!}e^{-\frac{i\gamma\pi k}{2\alpha}}\Gamma(\frac{k}{\alpha})
+\text{Re}\left(\frac{x^n}{\pi\alpha n!}\Gamma(\frac{n+1}{\alpha})e^{-\frac{i\gamma\pi(n+1)}{2\alpha}}R\right)\nn\\
&=\frac{-1}{\alpha\pi x}\text{Re}\sum_{k=1}^{n}\frac{(-x)^{k}}{(k-1)!}\exp\left(-i(\frac{k\pi(\gamma-\alpha)}{2\alpha} +\frac{\pi}{2})\right)\Gamma(\frac{k}{\alpha})+\text{Re}\left(\frac{x^n}{\pi\alpha n!}\Gamma(\frac{n+1}{\alpha})e^{-\frac{i\gamma\pi(n+1)}{2\alpha}}R\right)\nn\\
&:=\frac{1}{\pi x}\sum_{k=1}^{n}\frac{(-x)^{k}}{k!}\sin(\frac{k\pi(\gamma-\alpha)}{2\alpha})\Gamma(1+\frac{k}{\alpha})+R(x,\alpha),
\end{align}
where 
\beq\label{eqn:remainder}
\left| R(x,\alpha)\right|\leq \frac{x^n}{\pi\alpha n!}\Gamma(\frac{n+1}{\alpha}).
\eeq

It follows from Zolotarev's identity (see \cite[(7.10)]{Ko}) for $\alpha\in(\frac12,1)\cup (1,2)$,
$$
p^{(\alpha)}(x,\gamma)=x^{-1-\alpha}p^{(1/\alpha)}(x^{-\alpha}, \frac{\gamma+1}{\alpha}-1).
$$
Hence, from \eqref{eqn:HK} we have 
\begin{eqnarray*}
&&p^{(\alpha)}(x,\gamma)=x^{-1-\alpha}p^{(1/\alpha)}(x^{-\alpha}, \frac{\gamma+1}{\alpha}-1)\\
&=&x^{-1-\alpha}\left(\frac{1}{\pi x^{-\alpha}}\sum_{k=1}^{n}\frac{(-x^{-\alpha})^{k}}{k!}\sin(\frac{k\pi(\gamma-\alpha)}{2})\Gamma(1+k\alpha) +R(x^{-\alpha}, \frac{1}{\alpha})\right).
\end{eqnarray*}
Now let $\gamma=0$ and define $E(x,\alpha):=x^{-1-\alpha}R(x^{-\alpha}, \frac{1}{\alpha})$. From \eqref{eqn:remainder} we have 
$$
\left|E(x,\alpha)\right|\leq \frac{\alpha x^{-(n+1)\alpha-1}}{\pi n!}\Gamma(\alpha(n+1)).
$$
Hence, for all $\alpha\in (1,2)$ and $x\geq 1$ we have 
$$
\left|E(x,\alpha)\right|\leq \frac{2 x^{-n-2}}{\pi n!}\Gamma(2(n+1)).
$$

The fact $\frac{\Gamma(1+\alpha)\sin\frac{\pi\alpha}{2}}{\pi}=C\alpha$ can be proved using  a similar argument as in the proof of Lemma \ref{lemma:uniform ub}.
\qed

\begin{lemma}\label{lemma:Cauchy cancel}
Let $\alpha\in (1,2)$ and suppose that 
$Y^{(\alpha)}\in C_{k,k+1}$ for some $k\in \mathbb{N}$. 
There exists a function $\phi$, independent of $\alpha\in (1,2)$, such that 
$$
\Prob(\overline{Y}^{(\alpha)}_1>u, Y^{(\alpha)}_1\leq u)\leq \phi(u) \text{ with } \int_{0}^{\infty}\phi(u)du<\infty.
$$
\end{lemma}
\pf
Note that 
$$
\Prob(\overline{Y}^{(\alpha)}_{1}>u, Y_{1}^{(\alpha)}\leq u)=\Prob(\overline{Y}^{(\alpha)}_1>u)-\Prob(Y_{1}^{(\alpha)}>u).
$$
It follows from Lemmas \ref{lemma:uniform ub} and \ref{lemma:uniform lb} that
$$
\overline{p}^{(\alpha)}(x)\leq C\alpha x^{-1-\alpha} +
Ax^{-3} \text{ and } 
p^{(\alpha)}(x)\geq C\alpha x^{-1-\alpha} -\frac{12}{\pi}x^{-3},
$$
where $A$ is independent of $\alpha\in (1,2)$.

Hence, we have 
\begin{eqnarray*}
&&\Prob(\overline{Y}^{(\alpha)}_{1}>u, Y_{1}^{(\alpha)}\leq u)=\Prob(\overline{Y}^{(\alpha)}_1>u)-\Prob(Y_{1}^{(\alpha)}>u)\\
&\leq&\int_{u}^{\infty}\left(C\alpha x^{-1-\alpha} +
Ax^{-3}-C\alpha x^{-1-\alpha} +\frac{12}{\pi}x^{-3} \right)dx=\frac{A}{2}u^{-2}+\frac{6}{\pi}u^{-2}.
\end{eqnarray*}
Finally, we set 
$$
\phi(u):=1_{\{0<u\leq 1\}} + 1_{\{u>1\}}
(\frac{A}{2}+\frac{6}{\pi})u^{-2}.
$$
\qed

\begin{proposition}\label{prop:Cauchy ub}
Let $D$ be a bounded $C^{1,1}$ open set in $\R^{d}$, $d\geq 2$. Then, we have
$$
\limsup_{t\to 0}\frac{|D|-Q_{D}^{(1)}(t)}{t\ln(1/t)}\leq \frac{|\partial D|}{\pi}.
$$
\end{proposition}
\pf
From Lemma \ref{lemma:Cauchy ub} for any $\alpha\in(1,2)$ we have 
\beq\label{eqn:Cauchy ub main}
\limsup_{t\to 0}\frac{|D|-Q_{D}^{(1)}(t)}{t\ln(1/t)}\leq \frac{|\partial D|}{\pi} +
\frac{\int_{0}^{\infty}\Prob(\overline{Y}^{(\alpha)}_{1}>u, Y_{1}^{(\alpha)}<u )du }{\Gamma(1-\frac{1}{\alpha})}\cdot |\partial D|.
\eeq

Note that 
$$
\Gamma(x)=\int_{0}^{\infty}e^{-u}u^{x-1}du\geq \int_{0}^{1}e^{-u}u^{x-1}du\geq e^{-1}\int_{0}^{1}u^{x-1}du=\frac{1}{ex},
$$
and this implies 
\beq\label{eqn:Gamma}
\frac{1}{\Gamma(1-1/\alpha)}\leq e(1-\frac{1}{\alpha})=\frac{e(\alpha-1)}{\alpha} \text{ for all } \alpha>1.
\eeq
Note that from \eqref{eqn:ckl} for 
$Y^{(\alpha)}\in C_{k,k+1}$
we have $\alpha \in (1,2)$ and as $k\to \infty$, $\alpha\downarrow 1$.
Hence, it follows from Lemma \ref{lemma:Cauchy cancel} and the Lebesgue dominated convergence theorem that
$$
\lim_{Y_{t}^{(\alpha)}\in C_{k,k+1}, k\to \infty}\int_{0}^{\infty}\Prob(\overline{Y}^{(\alpha)}_{1}>u, Y_{1}^{(\alpha)}\leq u)du
=\int_{0}^{\infty}\Prob(\overline{Y}^{(1)}_{1}>u, Y_{1}^{(1)}\leq u)du.
$$
Note that the density $p^{(1)}(x)$ of $Y^{(1)}_{1}$ is given by 
$$
p^{(1)}(x)=\frac{1}{\pi(1+x^2)}
$$
and $\int_{u}^{\infty}\Prob(Y_{1}^{(1)}>u)du=\int_{u}^{\infty}\frac{1}{\pi (1+x^2)}dx=\frac{1}{\pi}(\frac{\pi}{2}-\arctan(u))=\frac{1}{\pi}\arctan(\frac{1}{u})$.
By the Taylor expansion of $\arctan x =\sum_{n=1}^{\infty}(-1)^{n+1}\frac{x^{2n-1}}{2n-1}$ for $|x|<1$ we have $\frac{1}{\pi}\arctan(1/u)\geq \frac{1}{\pi u} -\frac{1}{2\pi u^3}$ for all sufficiently large $u$.
Hence, it follows from \cite[(3.5)]{P20} that 
\begin{eqnarray*}
&&\Prob(\overline{Y}^{(1)}_{1}>u, Y_{1}^{(1)}\leq u ) =\Prob(\overline{Y}^{(1)}_{1}>u)-\Prob(Y_{1}^{(1)}> u)\\
&=&\Prob(\overline{Y}^{(1)}_{1}>u)-\frac{\arctan(1/u)}{\pi}\leq\Prob(\overline{Y}^{(1)}_{1}>u)-\frac{1}{\pi u} +\frac{1}{2\pi u^3}\leq \frac{4}{\pi^2}\frac{\ln u}{u^3} +\frac{1}{2\pi u^3}
\end{eqnarray*}
for all sufficiently large $u>0$, and this implies 
\beq\label{eqn:Cauchy finite}
\int_{0}^{\infty}\Prob(\overline{Y}^{(1)}_{1}>u, Y_{1}^{(1)}\leq u)du<\infty.
\eeq

Hence, if follows from \eqref{eqn:Gamma} and \eqref{eqn:Cauchy finite} that
$$
\lim_{Y^{(\alpha)}\in C_{k,k+1}, k\to \infty}
\frac{\int_{0}^{\infty}\Prob(\overline{Y}^{(\alpha)}_{1}>u, Y_{1}^{(\alpha)}<u )du }{\Gamma(1-\frac{1}{\alpha})}=0.
$$
Hence, from \eqref{eqn:Cauchy ub main} we reach the conclusion of the theorem.
\qed

\begin{theorem}\label{thm:main1}
Let $D$ be a bounded $C^{1,1}$ open set in $\R^{d}$, $d\geq 2$. Then, we have
$$
\lim_{t\to 0}\frac{|D|-Q_{D}^{(1)}(t)}{t\ln(1/t)}= \frac{|\partial D|}{\pi}.
$$
\end{theorem}
\pf
It follows immediately from \eqref{eqn:Cauchy lb} and Proposition \ref{prop:Cauchy ub}.
\qed

\vspace{0.2in}
\noindent
\textbf{Acknowledgments:} 
We thank the referee for helpful comments and suggestions.
The first-named author is grateful to Professor Kyeong Hun Kim (Korea University, South Korea) for his encouragement while this project was in progress. 
Part of the research for this project was done while the second-named author was visiting Jiangsu Normal University, where he was partially supported from the National Natural Science Foundation of China (11931004) and by the Priority Academic Program Development of Jiangsu Higher Education Institutions.
\bigskip  

\begin{remark}
\rm{After the paper was published in \textit{Electronic Journal of Probability} Vol. 27, paper no. 22, 1--19, (2022)}, the authors found a minor mistake in the proof of lower bound in Section \ref{section:alpha12}. 
In the proof of Lemma 3.1 there, it was stated that 
$$
\{\overline{Y}_{t}^{(\alpha)}>t\}=\{\tau_{H_{x}}^{(\alpha)}\leq t\}\subset \{\tau_{D}^{(\alpha)}\leq t\}. 
$$
However, the last inclusion is not true in general when $D$ is non-convex. We provide the corrected proof in Lemma \ref{lemma:lb}.
\end{remark}

\begin{singlespace}

\end{singlespace}
\end{doublespace}

\vskip 0.3truein

{\bf Hyunchul Park}

Department of Mathematics, State University of New York at New Paltz, NY 12561,
USA

E-mail: \texttt{parkh@newpaltz.edu}

\bigskip

{\bf Renming Song}

Department of Mathematics, University of Illinois, Urbana, IL 61801,
USA

E-mail: \texttt{rsong@illinois.edu}

\end{document}